\documentclass[12pt,a4paper]{article}

\usepackage[latin1]{inputenc}
\usepackage{color}

\usepackage{bussproofs}
\usepackage{graphicx}
\usepackage{graphics}
\usepackage{latexsym}
\usepackage{amssymb,amsmath}
\usepackage{euscript}
\usepackage{enumerate}
\usepackage{amsmath}
\usepackage{amsfonts}
\usepackage{multicol}
\usepackage{etex}
\usepackage[all]{xy}
\usepackage{color}
\usepackage{graphicx}
\usepackage{tikz}
\usepackage[latin1]{inputenc}
\usepackage[colorlinks=true, pdfstartview=FitV, linkcolor=blue, citecolor=blue, urlcolor=blue, pagebackref,]{hyperref}

\newtheorem{teo}{Theorem}[section]
\newtheorem{coro}[teo]{Corollary}
\newtheorem{lema}[teo]{Lemma}
\newtheorem{prop}[teo]{Proposition}

\newtheorem{defi}[teo]{Definition}

\newtheorem{rem}[teo]{Remark}
\newenvironment{dem}{\noindent {\em Proof:} }{\hfill $\square$ \bigskip}
\newenvironment{prueba}[1]{\noindent\bf Proof of #1. \rm}{$\quad \hfill \square$ \bigskip}
\numberwithin{equation}{section}

\newcommand{\R}{\mathbb{R}}

\newcommand{\N}{\mathbb{N}}
\renewcommand{\H}{\mathcal{H}}

\newcommand{\supp}{\mathop{\text{\rm supp}}}
\newcommand{\alphat}{\tilde{\alpha}}
\newcommand{\deltat}{\tilde{\delta}}
\newcommand{\menort}{\lesssim}
\newcommand{\mayort}{\gtrsim}

\def\XXint#1#2#3{{\setbox0=\hbox{$#1{#2#3}{\int}$}
    \vcenter{\hbox{$#2#3$}}\kern-.5\wd0}}

\usepackage{fancyhdr}

\fancypagestyle{encabezadotesis}{
\fancyhead{}
\fancyhead[RO]{\em \nouppercase{\rightmark}}
\fancyhead[LE]{\bf \nouppercase{\leftmark}}
\cfoot{\thepage}
}

\begin{document}
\title{On optimal parameters involved with two-weighted estimates of commutators of singular and fractional integral operators}
\author{Gladis Pradolini - Jorgelina Recchi}
\date{}
\maketitle

\renewcommand{\thefootnote}{\fnsymbol{footnote}}
\footnotetext{
The  authors are supported by CONICET and,
UNL and UNS, respectively.
}
\begin{abstract}
In this paper we prove two-weighted norm estimates for higher order commutator of singular integral and fractional type operators between weighted $L^p$ and certain spaces that include Lipschitz, BMO and Morrey spaces. We also give the optimal parameters involved with these results, where the optimality is understood in the sense that the parameters defining the corresponding spaces belong to certain region out of which the classes of weights are satisfied by trivial weights. We also exhibit pairs of non-trivial weights in the optimal region satisfying the conditions required.
\end{abstract}

\section{Introduction}
It is well known the significative contribution that represent the continuity properties of different operators from Harmonic Analysis in the study of the regularity properties of the solutions of certain partial differential equations. There is a vast evidence of this fact and in this direction, the commutators of operators with symbols functions in certain adequate spaces play an important role (see, for example, \cite{BraCeru}, \cite{BCe}, \cite{BraCe}, \cite{BCM}, \cite{ChFL}, \cite{ChFL2} and \cite{Rios}). Thus, their boundedness properties allow to derive regularity properties related with the solutions of such PDE's.

In \cite{HSV2} the authors proved one-weight boundedness results for the classical fractional integral operator $I_{\alpha}$, $0<\alpha<n$, including weighted estimates between $L^p$ spaces and certain generalizations of the $\mathcal{L}_{p,\lambda}$ spaces defined in \cite{Peetre}. Particularly, the Lipschitz spaces considered in that article are generalizations of some known integral version of the classical Lipschitz($\beta$) spaces, where the relation between $p$ and $\beta$ is the standard,  $\beta/n=\alpha/n-1/p$.
In \cite{Morvidone} similar problems were studied for the Hilbert transform and certain generalizations of the Lipschitz spaces defined in \cite{HSV2}.

On the other hand, in \cite{PradoCarolinae} a two-weighted problem for the boundedness of $I_{\alpha}$ in the spirit of \cite{HSV2} was studied. The parameters involved belong to a region out of which the weights are trivial, that is $v=0$ or $w=\infty$ a.e. and, in this sense, this is an optimal estimate. Similar results of this type were proved in \cite{PraReinf} for commutators of singular integral and fractional type operators.

In this paper we prove two-weighted norm estimates for singular integral and fractional type operators and their higher order commutators between weighted $L^p$ and certain spaces related to a parameter $\beta$, that include Lipschitz, BMO and Morrey spaces and that are wider than those considered in \cite{PraReinf}. Moreover the classes of weights are quite different from those given there, including local and global conditions. We also give the optimal parameters involved with these results, where the optimality is understood in the sense that the parameters $p$ and $\beta$ belong to certain region out of which the classes of weights are satisfied by trivial weights. Moreover, we exhibit concrete pairs of non-trivial weights in the optimal region satisfying the conditions required on the weights, where the boundedness results includes values of $\beta$ describing Lipschitz$(\beta)$, BMO and Morrey spaces, that is, $0<\beta<1$, $\beta=0$ and $\beta<0$, respectively. Our results extend those contained in \cite{PradoCarolinae} for the fractional integral operator (see also \cite{HSV2} for the one-weight case). We prove that a one-weight result can only holds whenever the relation between the parameters is standard.  We also study the relation between our classes of weights and those given in \cite{PraReinf}, which are natural extensions of the $A_1$-Muckenhoupt class in the one-weight estimates.

The paper is organized as follows. In section \S 2 we give the preliminaries and state the main results. In \S 3 we prove the optimality of the classes of weights and prove some other properties. Finally, in \S 4 we prove the main results.

\section{Preliminaries and main results}\label{preliminaries}

In this section we give  the definitions of  the  operators we shall be dealing with.
We shall consider singular integral operators of convolution type $T$ with kernel $K$, that is $T$ is bounded on $L^2(\R^n)$ and if $x\notin \supp f$
\begin{equation} \label{integral singular}
Tf(x)=\int_{\R^n}K(x-y)f(y)dy.
\end{equation}
The kernel $K$ is a measurable function defined away from $0$, satisfying certain smoothness condition to be described later. We shall also suppose that $K$  satisfies the  typical size condition  given by
$$
|K(x-y)|\leq \frac{C}{|x-y|^{n}},
$$
with will by called $S^*_0$.

Related with the singular integral operator $T$, we can formally define the commutator with symbol $b\in L^1_{{\rm loc}}(\R^n)$, by
$$
[b,T]f=b\,Tf-T(bf).
$$
The  commutator of order $m\in \N \cup \{0\}$ of $T$ is defined by
$$
T_b^0=T, \;\;\;\;\; T^m_b=[b,T^{m-1}_b].
$$
We shall also consider  fractional operators of convolution type $T_\alpha$, $0< \alpha<n$, defined by

\begin{equation}\label{operador fraccionario}
T_\alpha f(x)=\int_{\R^n}K_\alpha(x-y)f(y)dy,
\end{equation}
where the kernel $K_\alpha$ is not identically zero and verifies certain size and smoothness conditions.

Let $0<\delta< 1$. We say that a function $b$ belongs to the space $\Lambda(\delta)$ if there exists a positive constant $C$ such that, for every $x,y \in \R^n$
$$
|b(x)-b(y)|\leq  C |x-y|^\delta.
$$
The smallest of such constants will be denoted by $\|b\|_{\Lambda(\delta)}$. The space $\Lambda(\delta)$ is the well known Lipschitz space in the classical literature. We shall be dealing with commutators  with symbols belonging to this class of functions.



We say that $\Phi:[0,\infty)\to [0,\infty)$ is a Young function if it is increasing, convex and verifies $\Phi(0)=0$ and $\Phi(t)\to \infty$ when $t\to \infty$. The $\Phi$-Luxemburg average of a locally integrable function $f$ over a ball $B$ is defined by
$$
\|f\|_{\Phi,B}= \inf \left\{\lambda>0: \frac{1}{|B|}\int_B \Phi\left(\frac{|f(x)|}{\lambda}\right)\, dx \leq 1 \right\}.
$$

Given a Young function $\Phi$, the following H\"older's type inequality holds for every pair of measurable functions $f$, $g$
$$
\frac{1}{|B|}\int_B |f(x)g(x)|dx\leq 2 \|f\|_{\Phi,B}\|g\|_{\tilde{\Phi},B},
$$
where $\tilde{\Phi}$ is the complementary Young function of $\Phi$, defined by
$$
\tilde{\Phi}(t)=\sup_{s>0}\{st-\Phi(s)\}.
$$
It is well known that $t\leq\Phi^{-1}(t)\tilde{\Phi}^{-1}(t)\leq 2t$ for every $t>0$. Moreover, given $\Phi$, $\Psi$ and $\Theta$ Young functions verifying that $\Phi^{-1}(t)\Psi^{-1}(t)\menort \Theta^{-1}(t)$ for every $t>0$, the following generalization holds
$$
\|fg\|_{\Theta,B}\menort\|f\|_{\Phi,B}\|g\|_{\Psi,B}.
$$
For more information about Orlicz spaces see \cite{RR}.

We say that a  kernel $K_\alpha \in S_{\alpha,}^*$ with $0\leq \alpha<n$ , if there exists a positive constant $C$ such that,

$$
|K_\alpha(x)|\leq \frac{C}{|x|^{n-\alpha}}.
$$

We classify the operators $T$ and $T_{\alpha}$ into different types, according to the smoothness condition satisfied by   $K_\alpha$.
\subsection{Operators with Lipschitz regularity}

%

Let $0\leq \alpha<n$, we say that a kernel  $K_\alpha$  belongs to $K^{*}_{\alpha, \infty}$ if there exists a positive constant $C$ and $0<\eta\leq 1$ such that
$$
|K_\alpha(x-y)-K_\alpha(x'-y)|+|K_\alpha(y-x)-K_\alpha(y-x')|\leq C \frac{|x-x'|^\eta}{|x-y|^{n-\alpha+\eta}},
$$
whenever $|x-y|\geq 2|x-x'|$. We denote $K_0=K$, the kernel in \eqref{integral singular} and $T_0=T$ the singular integral operator. 

It is easy to check that the fractional integral operator $I_\alpha$, with kernel $K_\alpha(x)=|x|^{\alpha-n}$, satisfies conditions $S^*_{\alpha}$
and $K^{*}_{\alpha, \infty}$ for $0<\alpha<n$. On the other hand, it it well known that the Hilbert transform satisfies $S^*_{0}$ and $K^{*}_{0, \infty}$.

Related with the fractional type integral operators $T_\alpha$, we can formally define the higher order commutators with symbol $b\in L^1_{{\rm loc}}(\R^n)$, by
$$
T_{\alpha,b}^m f(x)=\int_{\R^n}(b(x)-b(y))^m K_\alpha (x-y)f(y)\,dy,
$$
where $m\in \N \cup \{0\}$ is the order of the commutator. Clearly, $T^0_{\alpha,b}=T_\alpha$.

As we have said, we are  interested in studying the boundedness properties of the commutators $T_{\alpha,b}^m$ with symbol $b\in \Lambda(\delta)$,  from weighted Lebesgue spaces into certain weighted version ${\mathcal L}_{p,\lambda}$ spaces introduced in \cite{Peetre}. For $\beta \in \R $ and a weight $w$, these spaces are denoted by $\mathcal{L}_w(\beta)$ and collect the functions $f\in L^{1}_{{\rm loc}}(\R^n)$ that satisfy
$$
\frac{1}{w(B)|B|^{\beta}}\int_B |f(x)-m_B(f)|\,dx\,\leq\,C.
$$
for some positive constant $C$.
When $\beta=0$, $\mathcal{L}_w(0)$ is a weighted version of the bounded  mean oscilation  space introduced by Muckenhoupt and Wheeden in \cite{MW2}. Moreover,  $\mathcal{L}_1(\beta)$ gives the known Lipschitz integral space for $\beta$ in the range $0<\beta<1/n$ and the Morrey space, for $-1<\beta<0$. This class of functions was introduced in \cite{HSV2}.

In \cite{PraReinf} the authors prove two weighted boundedness results for commutators of a great variety of operators between Lebesgue and Lipschitz spaces $\mathbb{L}_w(\beta)$. These spaces collect the functions $f\in L^{1}_{{\rm loc}}(\R^n)$ that satisfy
$$
\frac{\|(1/w)\chi_B\|_{\infty}}{|B|^{1+\beta}}\int_B |f(x)-m_B(f)|\,dx\,\leq\,C.
$$
for some positive constant $C$. It is easy to check that, for a general weight $w$, $\mathbb{L}_w(\beta)\subset \mathcal{L}_w$ and, if $w$ belongs to the $A_1$-Muckenhoupt class then both spaces coincide.

In the sequel we denote $\alphat=m\delta+\alpha$. We say that $A\menort B$ if there exists a positive constant $c$ such that $A\leq c\,B$.

Related to the spaces $\mathcal{L}_w(\beta)$, we introduce the following class of weights.



%

\begin{defi} \label{definicion de la clase}
Let  $0\leq\alpha<n$, $0\leq\delta\leq 1$ and $1 <r\leq \infty$. Let $m\in \N \cup \{0\}$, and  $\deltat\leq\delta$. We say that a pair of weights $(w,v)$ belongs to $\H (r,\alphat,\deltat)$, if the inequality
\begin{equation}\label{condicion del peso}
|B|^{\frac{\delta-\deltat}{n}}\left(\int_{\R^n}\frac{v^{r'}(y)}{\left(|B|^{1/n}+|x_B-y|\right)^{r'(n-\alphat+\delta)}}\,dy\right)^{1/{r'}} \menort \frac{w(B)}{|B|}
\end{equation}
holds for every ball $B\subset \R^n$, where $x_B$ is the center of $B$. In the case $r=1$ we say that $(w,v)$ belongs to $\H(1,\alphat,\deltat)$ if the inequality
\begin{equation}\label{condicion del peso para r=1}
|B|^{\frac{\delta-\deltat}{n}}\left\|\frac{v(.)}{\left(|B|^{1/n}+|x_B-.|\right)^{n-\alphat+\delta}}\right\|_{\infty} \menort \frac{w(B)}{|B|}
\end{equation}
holds for every ball $B\subset \R^n$, where $x_B$ is the center of $B$.

\end{defi}
When $0<\alpha<n$, $m=0$ and $\delta=1$, the class of pair of weights in Definition $\ref{definicion de la clase}$ was introduced in \cite{PradoCarolinae}. If in adittion, $w=v$ and $\deltat=\alphat-n/r$, then the class $\H (r,\alphat,\deltat)$ was  defined in \cite{HSV2}. When $n=1$, $\alpha=0$, $\delta=1$, $m=0$ and $w=v$ the class $\H(\infty,0,0)$ is the class $B_2$ in \cite{MW2}.
\begin{rem}\label{clases contenidas}
Since
$$
\|(1/w)\chi_B\|_{\infty}= \frac{1}{\inf_{x\in B}w}\geq \frac{|B|}{w(B)}
$$
then the classes $\mathbb{H}(r,\alphat,\deltat)$ defined in \cite{PraReinf} are contained in the classes $\H(r,\alphat,\deltat)$. We shall prove later that this inclusion is strict.
\end{rem}

\begin{rem}
We say $w\in \H(r,\alphat,\deltat)$ if $w=v$ in (\ref{condicion del peso}) and (\ref{condicion del peso para r=1}).
\end{rem}
In the one-weight case, we obtain the following lemma.

\begin{lema} \label{un solo peso solo en el borde}
Let  $0\leq\alpha<n$, $0<\delta<1$ and $1\leq r\leq \infty$. Let $\deltat\leq\min\{\delta,\alphat-n/r\}$. If $w\in \H (r,\alphat,\deltat)$, then $\deltat=\alphat-n/r$.
\end{lema}

\begin{dem}
Let $1\leq r\leq \infty$ (if $r=1$ we understand $\|.\|_{\infty}$ instead $\|.\|_{r'}$). Since $w\in \H (r,\alphat,\deltat)$, we have
\begin{align*}
\frac{w(B)}{|B|}&\gtrsim |B|^{\frac{\delta-\deltat}{n}}\left(\int_{B}\frac{w^{r'}(y)}{\left(|B|^{1/n}+|x_B-y|\right)^{r'(n-\alphat+\delta)}}\,dy\right)^{1/{r'}}\\
& \mayort |B|^{\frac{\delta-\deltat}{n}-(\frac{1}{r}-\frac{\alphat}{n}+\frac{\delta}{n})}\left(\frac{1}{|B|}\int_{B}w^{r'}(y)\,dy\right)^{1/{r'}}\\
& \mayort |B|^{-\frac{\deltat}{n}-\frac{1}{r}+\frac{\alphat}{n}}\frac{w(B)}{|B|}.
\end{align*}
Then, this inequality holds if $\deltat=\alphat-n/r$.
\end{dem}

\space

We are now in a position to state our main results. We first state the theorem for singular integral operators with the corresponding weights belonging to $\H(r,m\delta,\deltat)$, that is $\alpha=0$ in Definition \ref{definicion de la clase}. 

We say that a weight $w$ belongs to the reverse H\"{o}lder class $RH(s)$ if there exists a positive constant $C$ such that
$$
\left(\frac{1}{|B|}\int_B w^{s}(x)dx\right)^{1/{s}}\leq C \frac{w(B)}{|B|}.
$$

\begin{teo} \label{teo para integrales singulares}
Let $0<\delta<\min \{\eta,n/m\}$ and $1\leq r \leq \infty$. Let $\deltat\leq \min\{\delta, m\delta-n/r\}$ and $b\in \Lambda(\delta)$. If $(w,v) \in \H(r,m\delta,\deltat)$,  $v^{r'} \in RH(\frac{p(r-1)}{r-p})$ for some $p\in (1,r)$  and $K \in S^*_0$, then
$$
\|T^m_b f\|_{\mathcal{L}_w(\deltat/n)} \menort\,\|b\|^m_{\Lambda(\delta)}\left\| f/v \right\|_{L^r(\R^n)},
$$
holds for every $f$ such that $f/v \in L^r(\R^n)$.
\end{teo}
From the theorem above and Lemma \ref{un solo peso solo en el borde} we obtain the following corollary.
\begin{coro} \label{teo para integrales singulares con un solo peso} Let $0<\delta<\min \{\eta,n/m\}$ and $1\leq r <\infty$. Let $\deltat= m\delta-n/r$ and $b\in \Lambda(\delta)$. If $w \in \H(r,m\delta,\deltat)$, and $K \in S^*_0$, then
$$
\|T^m_b f\|_{\mathcal{L}_w(\deltat/n)} \menort\,\|b\|^m_{\Lambda(\delta)}\left\| f/w \right\|_{L^r(\R^n)},
$$
holds for every $f$ such that $f/w \in L^r(\R^n)$.
\end{coro}

For the Hilbert transform, $m=0$ and $r=\infty$, this corollary was proved in \cite{MW2}.

\medskip

We now state the main results for the boundedner Fractional integral operators, that is $0<\alpha<n$.
\begin{teo}\label{teorema A}
Let $0<\alpha<n$, $0<\delta<\min\{\eta, (n-\alpha)/m\}$ and $1\leq r\leq \infty$. Let $\deltat\leq \min\{\delta,\alphat-n/r\}$ and $b\in \Lambda(\delta)$. If  $(w,v)\in \H(r,\alphat, \deltat)$ and $K_\alpha \in S^*_\alpha$, then
$$
\|T^m_{\alpha,b}f\|_{\mathcal{L}_w(\deltat/n)} \menort\,\|b\|^m_{\Lambda(\delta)}\left\| f/v \right\|_{L^r(\R^n)},
$$
holds for every $f$ such that $f/v \in L^r(\R^n)$.
\end{teo}
From Theorem \ref{teorema A} and Lemma \ref{un solo peso solo en el borde}, we obtain the following result.

\begin{coro}
Let $0<\alpha<n$, $0<\delta<\min\{\eta, (n-\alpha)/m\}$ and $1\leq r\leq \infty$. Let  $\deltat=\alphat-n/r$ and $b\in \Lambda(\delta)$. If  $w\in \H(r,\alphat,\deltat)$ and $K_\alpha \in S^*_\alpha$, then
$$
\|T^m_{\alpha,b}f\|_{\mathcal{L}_w(\deltat/n)} \menort\,\|b\|^m_{\Lambda(\delta)}\left\| f/w \right\|_{L^r(\R^n)},
$$
holds for every $f$ such that $f/w \in L^r(\R^n)$.
\end{coro}


%
%

\section{Properties of the classes of weights }

We  give some properties of the  classes of weights $\H (r,\alphat,\deltat)$ given in definition \ref{definicion de la clase} . Recall that $\deltat\leq \min\{ \alphat-n/r,\delta\}$ and $\alphat=m\delta+\alpha$, where $0\leq \alpha< n$; $1\leq r\leq \infty$.

In this section we shall proof that the range of the parameters involved in the classes $\H(r,\alphat,\deltat)$  lie in the  shaded region of the Figure 1. 

\begin{figure}[h]
	\begin{tikzpicture}[xscale=0.5,yscale=0.8]	
	\filldraw[fill=gray!30,draw=white] (0,-4.5)--(0,3)--(2.5,3)--(5,-2)--(5,-4.5)--cycle;
	\draw (0,-4.5)--(0,3)--(2.5,3)--(5,-2)--(5,-4.5);
	\draw[->] (-1,0)--(7,0) node[below] {$1/r$};
	\draw[->] (0,-4.5)--(0,4) node[left] {$\tilde{\delta}$};
	\draw[fill=white] (2.5,3) circle [radius=0.08];
	\node[left] at (0,3) {$\delta$};
	\draw[dashed] (5,-2)--(0,-2) node[left] {$\tilde{\alpha}-n$};
	\draw (5,-0.1)--(5,0.1);
	\node[below] at (5,0) {$1$};
	\node at (4,4) {$\tilde{\alpha} > 1$};
	\draw [<-] (3.5,1.5) -- ++(0.2,6pt) -- ++(0.2,-6pt) -- ++(0.2,6pt) node[right] {$\tilde{\delta} = \tilde{\alpha} - n/r$};
	\end{tikzpicture}
	\begin{tikzpicture}	[xscale=0.5,yscale=0.8]	\hspace{-1cm}
	\filldraw[fill=gray!30,draw=white] (0,-4.5)--(0,3)--(5,-2)--(5,-4.5)--cycle;
	\draw (0,-4.5)--(0,3)--(5,-2)--(5,-4.5); 
	\draw[->] (-1,0)--(7,0) node[below] {$1/r$}; 
	\draw[->] (0,-4.5)--(0,4) node[left] {$\tilde{\delta}$}; 
	\draw[fill=white] (0,3) circle [radius=0.08]; 
	\node[left] at (0,3) {$\delta$}; 
	\draw[dashed] (5,-2)--(0,-2) node[left] {$\tilde{\alpha}-n$}; 
	\draw (5,-0.1)--(5,0.1); 
	\node[below] at (5,0) {$1$}; 
	\node at (4,4) {$\tilde{\alpha} = 1$};
	\hspace{-0.8cm}	\draw [<-] (3.5,1.5) -- ++(0.2,6pt) -- ++(0.2,-6pt) -- ++(0.2,6pt) node[right] {$\tilde{\delta} = 1 - n/r$};
	\end{tikzpicture}\begin{tikzpicture}[xscale=0.5,yscale=0.8]	 \hspace{-2cm}
	\filldraw[fill=gray!30,draw=white] (0,-4.5)--(0,2.5)--(5,-2)--(5,-4.5)--cycle;
	\draw (0,-4.5)--(0,2.5)--(5,-2)--(5,-4.5);
	\draw[->] (-1,0)--(7,0) node[below] {$1/r$};
	\draw[->] (0,-4.5)--(0,4) node[left] {$\tilde{\delta}$};
	\draw (-0.1,3)--(0.1,3);
	\node[left] at (0,3) {$\delta$};
	\draw[dashed] (5,-2)--(0,-2) node[left] {$\tilde{\alpha}-n$};
	\draw (5,-0.1)--(5,0.1);
	\node[below] at (5,0) {$1$};
	\node at (4,4) {$\tilde{\alpha} < 1$};
	\hspace{-1cm}\draw [<-] (3.5,1.5) -- ++(0.2,6pt) -- ++(0.2,-6pt) -- ++(0.2,6pt) node[right] {$\tilde{\delta} = \tilde{\alpha} - n/r$};
	\end{tikzpicture}
	\caption{Permissible range of the parameters $r$ and $\deltat$ for different values of $\alphat$.}
	\label{h}
\end{figure}
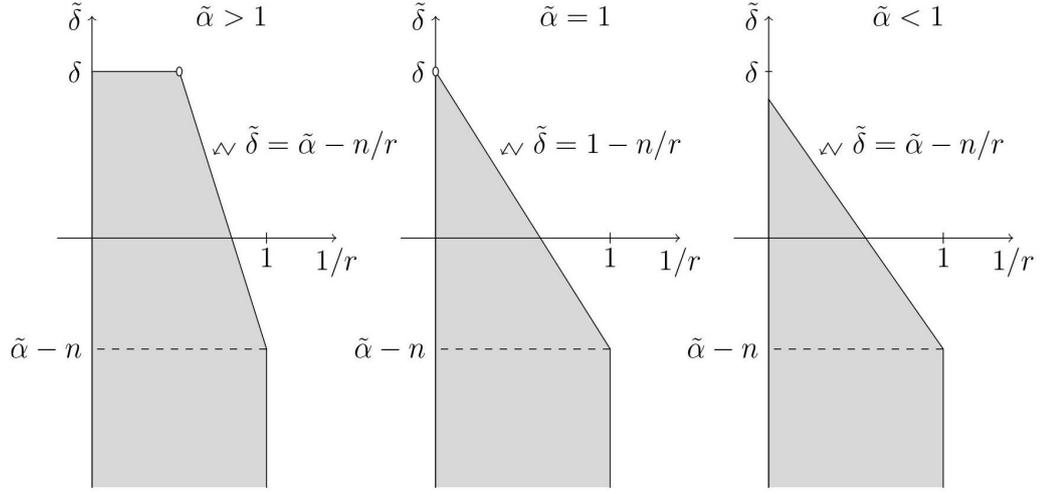
\begin{lema}\label{lema del doble v}
Let  $0\leq\alpha<n$, $0<\delta<\min\{\eta, (n-\alpha)/m\}$ and $1\leq r\leq \infty$. If $(w,v) \in\H (r,\alphat,\deltat)$. Then
$$
\|v\chi_{2B}\|_{r'}\menort |B|^{\frac{\deltat-\alphat}{n}}w(B).
$$

\end{lema}
\begin{dem}
Since $(w,v)\in \H(r, \alphat,\deltat)$, we know that
$$
|B|^{(\delta-\deltat)/n}\left(\int_{\R^n}\frac{v^{r'}(y)}{(|B|^{1/n}+|x_B-y|)^{r'(n-\alphat+\delta)}}dy\right)^{1/{r'}}\leq C\frac{w(B)}{|B|}
$$
for every ball $B\subset\R^n$, where $x_B$ is the center of $B$. Then, we have
\begin{align*}
v^{r'}(2B)& =\frac{|B|^{(n-\alphat+\delta)r'/n}}{|B|^{(n-\alphat+\delta)r'/n}}\int_{2B}v^{r'}(y)dy\\
& \menort |B|^{(n-\alphat+\delta)r'/n}\int_{2B}\frac{v^{r'}(y)}{(|B|^{1/n}+|x_B-y|)^{r'(n-\alphat+\delta)}}dy\\
& \menort |B|^{(n-\alphat+\delta)r'/n}\int_{\R^n}\frac{v^{r'}(y)}{(|B|^{1/n}+|x_B-y|)^{r'(n-\alphat+\delta)}}dy\\
& \menort |B|^{(n-\alphat+\delta)r'/n}\left(\frac{w(B)}{|B|^{\frac{\delta-\deltat}{n}+1}}\right)^{r'}\\
& \menort |B|^{(\deltat-\alphat)r'/{n}}\left(w(B)\right)^{r'}.
\end{align*}
\end{dem}

As a consequence of the lemma above and Lemma \ref{un solo peso solo en el borde} we obtain the following result.
\begin{coro} Let  $0\leq\alpha<n$, $0<\delta<\min\{\eta, (n-\alpha)/m\}$ and $1\leq r\leq \infty$. If $w \in\H (r,\alphat,\deltat)$. Then $w^{r'}$ satisfies a doubling condition.
\end{coro}
When $m=0$ this corollary was proved in \cite{HSV2}.\\

As a consequence of  Lemmas  \ref{un solo peso solo en el borde} and \ref{lema del doble v}, we get the following result.

\begin{coro}\label{lema 3.8 de pola}
Let  $0\leq \alpha<n$, $0<\delta<\min\{\eta, (n-\alpha)/m\}$ and  $1\leq r\leq \infty$. If $w$ be a weight in $\H (r,\alphat,\deltat)$, then $w$ belongs to $RH(r')$.
\end{coro}
The next lemma proves the equivalence between the class $\H(r,\alphat,\deltat)$ with a local and global conditions. The proof is straightforward and we omit it.
\begin{lema}\label{lema condicion de peso es equivalente a la parte local mas la global}
Let  $0\leq \alpha<n$, $0<\delta<\min\{\eta, (n-\alpha)/m\}$ and  $1\leq r\leq \infty$. The condition $\H (r,\alphat,\deltat)$ is equivalent to  the  following two inequalities
\begin{equation}\label{parte local}
|B|^{\frac{\alphat-\deltat}{n}-\frac{1}{r}}\left(\frac{1}{|B|}\int_B v^{r'}(y)dy\right)^{1/{r'}}\menort \,\frac{w(B)}{|B|}
\end{equation}
and
\begin{equation}\label{parte global}
\displaystyle|B|^{\frac{\delta-\deltat}{n}}\left(\int_{\R^n-B}\frac{v^{r'}(y)}{|x_B-y|^{r'(n-\alphat+\delta)}}dy\right)^{1/{r'}}\menort \frac{w(B)}{|B|}
\end{equation}
hold simultaneously for every ball $B\subset \R^n$, where $x_B$ is the center of $B$.
\end{lema}

It is important to note that both condition (\ref{parte local}) and (\ref{parte global}) cannot be reduced to (\ref{parte global}) as in \cite{HSV2} for the one-weighted case. However, under certain additional hypothesis on $v$ then $\H (r,\alphat,\deltat)$ is condition  (\ref{parte global}). This fact established in the following lemma.

\begin{lema} \label{global implica local}
Let  $0\leq \alpha<n$, $0<\delta<\min\{\eta, (n-\alpha)/m\}$ and  $1\leq r\leq \infty$. Let $v$ be a weight such that $v^{r'}$ satisfies a doubling condition. Then, the global condition (\ref{parte global}) implies the local condition (\ref{parte local}).
\end{lema}
\begin{dem}
Since $v^{r'}$ is a doubling weight, we get
\begin{align*}
\left(\frac{1}{|B|}\int_B v^{r'}(y)dy\right)^{1/{r'}}& \menort \left(\frac{1}{|B|}\int_{2B-B} v^{r'}(y)dy\right)^{1/{r'}}\\
&= \frac{|B|^{(n-\alphat+\delta)/n-1/r'}}{|B|^{(n-\alphat+\delta)/n}}\left(\int_{2B-B} v^{r'}(y)dy\right)^{1/{r'}}\\
&\menort |B|^{(n-\alphat+\delta)/n-1/r'}\left(\int_{2B-B} \frac{v^{r'}(y)}{|x_B-y|^{r'(n-\alphat+\delta)}}dy\right)^{1/{r'}}.
\end{align*}
Then, using the global condition (\ref{parte global}), we have that
$$
\left(\frac{1}{|B|}\int_B v^{r'}(y)dy\right)^{1/{r'}}\menort |B|^{(\deltat-\alphat)/n+1/r} \frac{w(B)}{|B|},
$$
which is (\ref{parte local}).
\end{dem}

Even though $v^{r'}$ satisfies a doubling condition, both inequalities (\ref{parte local}) and (\ref{parte global}) are not equivalent. This fact is showed in Lemma \ref{lema pesos locales y no globales}. We first give well known estimates in order to prove it.
%
%
%

\begin{lema}\label{lema 3.5 de Gladis}
Let $B=B(x_B,R)\subset \R^n$ and $\alpha >-n$. They are true the following statements,
\begin{enumerate}
\item If $|x_B|\leq R$,  $\int_B |x|^{\alpha}dx\approx R^{\alpha+n}$
\item If $|x_B|>R$,  $\int_B|x|^{\alpha}dx\approx|x_B|^{\alpha}R^n$.
\end{enumerate}
\end{lema}

\begin{lema}\label{lema pesos locales y no globales}
 Let  $0\leq \alpha<n$, $0<\delta<\min\{\eta, (n-\alpha)/m\}$ and  $1\leq r\leq \infty$. There exist nontrivial pairs of weights $(w,v)$ that satisfy the local condition (\ref{parte local})  but not the global condition (\ref{parte global}) for $\deltat$ in the range
$$
\deltat\leq \min\{\delta, \alphat-n/r\},
$$
excluding the case $\deltat=\delta$ when $\alphat-n/r=\delta$.
\end{lema}
\begin{dem} Let us first consider $\deltat=\delta<\alphat-n/r$. Let $w=1$ and $v(x)=|x|^{n/r-\alphat+\delta}$, we shall proof that $(w,v)$ satisfies (\ref{parte local}) but not (\ref{parte global}). Indeed,
let $B=B(x_B,R)$. By Lemma \ref{lema 3.5 de Gladis}, if $|x_B|\leq R$  we get
$$
\frac{|B|^{(\alphat-\delta)/n}}{w(B)}\left(\int_B v^{r'}(y) dy\right)^{1/{r'}}\menort R^{\alphat-\delta-n}R^{n/r-\alphat+\delta+n/{r'}}\menort C,
$$
and if $|x_B|\geq R$ we have
$$
\frac{|B|^{(\alphat-\delta)/n}}{w(B)}\left(\int_B v^{r'}(y) dy\right)^{1/{r'}}\menort R^{\alphat-\delta-n}|x_B|^{n/r-\alphat+\delta}R^{n/{r'}}\menort C.
$$
On the other hand, if we now take $B=B(0,R)$, we get
\begin{align*}
\frac{|B|}{w(B)}\left(\int_{\R^n\setminus B} \frac{v^{r'}(y)}{|y|^{(n-\alphat+\delta)r'}} dy\right)^{1/{r'}}& \mayort  \left(\int_{\{|y|>R\}}\frac{|y|^{(n/r-\alphat+\delta)r'}}{|y|^{(n-\alphat+\delta)r'}} dy \right)^{1/{r'}}\\
&\mayort \left(\int_{\{|y|>R\}}\frac{1}{|y|^{n}} dy \right)^{1/{r'}}
\end{align*}
and the last integral is infinite. Thus, $(w,v)$ does not satisfy (\ref{parte global}).\\
Similar estimates can be obtained for the case $\deltat <\delta\leq \alphat-n/r $ by considering $(|x|^{\alphat-\deltat-n/r},1)$.  For the case $\deltat \leq \alphat-n/r\leq \delta$ the same is true for $(|x|^\beta,|x|^\theta)$, with $\theta>n/r-\alphat+\delta$ and $\beta=\theta+\alphat-\deltat-n/r$.
\end{dem}

\begin{prop}
Let $0\leq \alpha<n$, $0< \delta <\min\{\eta, (n-\alpha)/m\}$ and  $1\leq r\leq \infty$. Then
\begin{enumerate}[$i)$]
\item If $\deltat >\delta$ or $\deltat>\alphat-n/r$, then $(w,v)\in \H(r,\alphat,\deltat)$ if and only if $v=0$ almost everywhere in $\R^n$.
\item   If $\deltat=\alphat-n/r=\delta$  then the same conclusion as in $i)$ holds.
\end{enumerate}
\end{prop}
\begin{dem}
Let us first see $i)$ and let $\deltat >\delta$. Let $B=B(x,R)$ where $x$ is a Lebesgue point of $w$. Suposse that $r$ is finite, since  $(w,v)\in \H(r,\alphat,\deltat)$ we get
$$
\left(\int_{\R^n}\frac{v^{r'}(y)}{(|B|^{1/n}+|x-y|)^{r'(n-\alphat+\delta)}}dy\right)^{1/{r'}}\menort \frac{w(B)}{|B|}|B|^{(\deltat-\delta)/n}.
$$
From inequality above, by letting  $R \to 0$,  we obtain that
$$
\left(\int_{\R^n}\frac{v^{r'}(y)}{(|B|^{1/n}+|x-y|)^{r'(n-\alphat+\delta)}}dy\right)^{1/{r'}}=0,
$$
and so $v=0$ a.e. $x\in \R^n$.\\
Now, if  $\deltat>\alphat-n/r$, since $(w,v)\in \H(r,\alphat,\deltat)$ and by Lemma \ref{lema condicion de peso es equivalente a la parte local mas la global}, we have
$$
\left(\frac{v^{r'}(B)}{|B|}\right)^{1/r'}\menort \frac{w(B)}{|B|}|B|^{1/r-(\alphat-\deltat)/n}.
$$
If we choose $B(x,R)$ as before, since $x$ is a Lebesgue point of $w$, we get
$$
\lim_{R\to 0} \frac{w(B)}{|B|}|B|^{1/r-(\alphat-\deltat)/n}=0,
$$
from which it follows that
$$
\limsup_{R\to 0} \frac{v^{r'}(B(x,R))}{|B(x,R)|}=0.
$$

Clearly we can get the same conclusion for a.e. $x \in \R^n$. By standard arguments we can deduced that $v(x)=0$ in a.e. $x\in \R^n$.
If $r=\infty$  we have  to consider $1/r=0$ and $r'=1$ in the previous proof.

We now proceed with the proof of $ii)$, since $\deltat=\delta=\alphat-n/r$ we are going to see that   $(w,v)\in \H(r,\alphat,\deltat)$, with $r=n/{(\alphat-\delta)}$ if $\alphat >\delta$ and $r=\infty$ if $\alphat=\delta$, if $v(x) = 0$ in a.e..\\
Let $B=B(x_0, R)\subset \R^{n}$, since  $(w,v)\in \H(r,\alphat,\deltat)$  we get
$$
\left(\int_{\R^n}\frac{v^{r'}(y)}{(|B|^{1/n}+|x_0-y|)^{r'(n-\alphat+\delta)}}dy\right)^{1/{r'}}\menort \frac{w(B)}{|B|}.
$$
Since $n-\alphat+\delta=n/r'$, we have

\begin{equation}\label{A}
\left(\int_{\R^n}\frac{v^{r'}(y)}{(|B|^{1/n}+|x_0-y|)^{n}}dy\right)^{1/{r'}}\menort\frac{w(B)}{|B|}.
\end{equation}
We  now proceed as in the proof of Theorem $5.6$ in \cite{Pra} in orden to obtain that $v(x)=0$ a.e. $x\in \R^n$.
\end{dem}

\begin{rem}Let $0\leq \alpha<n$, $0<\delta<\min\{\eta, (n-\alpha/m)\}$ and $n/{\alphat}<r<n/{(\alphat-\delta)}$.  Let
$$
(2(\alphat-n/r)-\delta)^{+}\leq \theta \leq \alphat-n/r,
$$
$$
\alphat-n/r-\theta <\deltat <\min\{\alphat-n/r,n/r-\alphat+\delta\},
$$
we now exhibit a pair of weights $(w,v)$ such that $(w,v) \in \H(r,\alphat,\deltat)$ but $(w,v) \notin \mathbb{H}(r,\alphat,\deltat)$.
\end{rem}
Let
$w(x)=|x|^{\theta}\chi_{\{|x|\leq 1\}}+|x|^{\theta+\deltat}\chi_{\{|x| >1\}}$ and $v(x)=|x|^{\deltat}$. It is easy to check that $(w,v)$ does not belong to $\mathbb{H}(r,\alphat,\deltat)$. However, we shall see that $(w,v)\in \H(r,\alphat,\deltat)$. Since $v^{r'}$ is a doubling weight, by Lemma \ref{global implica local}, we shall only proof (\ref{parte global}).\\
Let $B=B(x_B,R)$ and $B_i=2^{i}B$. If $|x_B|\leq R$, by Lemma \ref{lema 3.5 de Gladis} we obtain that
\begin{align*}
\frac{|B|^{1+(\delta-\deltat)/n}}{w(B)}&\left(\int_{\R^n\setminus B}\frac{v^{r'}(y)}{|x_B-y|^{(n-\alphat+\delta)r'}}dy\right)^{1/{r'}}\\
&\menort \frac{R^{\alphat-\deltat}}{w(B)}\sum_{i=1}^{\infty}2^{-i(n-\alphat+\delta)}\left(v^{r'}(B_i)\right)^{1/{r'}}\\
&\menort \frac{R^{\alphat+n/{r'}}}{w(B)}\sum_{i=1}^{\infty}2^{-i(n/r-\deltat-\alphat+\delta)}\\
&\menort \frac{R^{\alphat+n/{r'}}}{w(B)}.
\end{align*}
Thus, since $w(B)\mayort \max\{R^{\theta+n},R^{\theta+\deltat+n}\}$ we obtain that (\ref{parte global}) holds for this case.\\
Now, if $|x_B|>R$, then there exists $N_1$ such that $2^{N_1}R\leq |x_B|\leq 2^{N_1+1}R$.
\begin{align*}
\frac{|B|^{1+(\delta-\deltat)/n}}{w(B)}&\left(\int_{\R^n\setminus B}\frac{v^{r'}(y)}{|x_B-y|^{(n-\alphat+\delta)r'}}dy\right)^{1/{r'}}\\
&\menort \frac{R^{\alphat-\deltat}}{w(B)}\sum_{i=1}^{\infty}2^{-i(n-\alphat+\delta)}\left(v^{r'}(B_i)\right)^{1/{r'}}\\
&= \frac{R^{\alphat-\deltat}}{w(B)}\sum_{i=1}^{N_1}2^{-i(n-\alphat+\delta)}\left(v^{r'}(B_i)\right)^{1/{r'}}+\frac{R^{\alphat-\deltat}}{w(B)}\sum_{i=N_1+1}^{\infty}2^{-i(n-\alphat+\delta)}\left(v^{r'}(B_i)\right)^{1/{r'}}\\
&=S_1+S_2.
\end{align*}
Let us first estimate $S_1$. Since $i\leq N_1$, $\;n/r-\alphat+\delta>0$ and  $w(B)\mayort \max \{|x_B|^{\alphat}R^{n},|x_B|^{\theta+\deltat}R^{n}\}$ we have
$$
S_1\menort \frac{R^{\alphat-\deltat+n/{r'}}}{w(B)}|x_B|^{\deltat}\menort C.
$$
In order to estimate $S_2$, we first observe that
$$
S_2\menort \frac{R^{\alphat+n/{r'}}}{w(B)}
$$
and then we proceed as in the estimate of $S_1$ to obtain that $S_2\menort C$.

\begin{teo}\label{teo 4.17 de tesis gladis}
Let $0\leq \alpha<n$ and  $0< \delta <\min\{\eta, (n-\alpha)/m\}$. There exist pairs of weights with $v$ not identically equal to zero, that verify the condition $\H(r,\alphat,\deltat)$ in the range of $r$ and $\deltat$ given by
$$
1\leq r\leq \infty \;\;\,\; \text{and} \;\;\,\; \deltat \leq \min\{\delta, \alphat-n/r\}
$$
excluding the case $\deltat=\delta$ when $\alphat-n/r=\delta$.
\end{teo}
\begin{dem}
By Remark \ref{clases contenidas}, the pair of weights given in \cite{PraReinf} belong to $\H(r,\alphat,\deltat)$, for $1\leq r\leq \infty$ and $\alphat-n\leq \deltat \leq \min\{\delta,\alphat-n/r\}$ excluding the case $\deltat=\delta$ when $\alphat-n/r=\delta$.

So, we shall exhibit  examples of pairs of weights for the case $\deltat<\alphat-n$. We first  consider $1 < r\leq \infty$. We divide the range $\deltat<\alphat-n$ in two regions:
\begin{enumerate}[$i)$]
\item $\alphat-n-k\delta<\deltat\leq \min \{\alphat-n/r-k\delta, \alphat-n-(k-1)\delta\}$, $k\in \N$. (see figure \ref{FIG: pesos en los dientes})
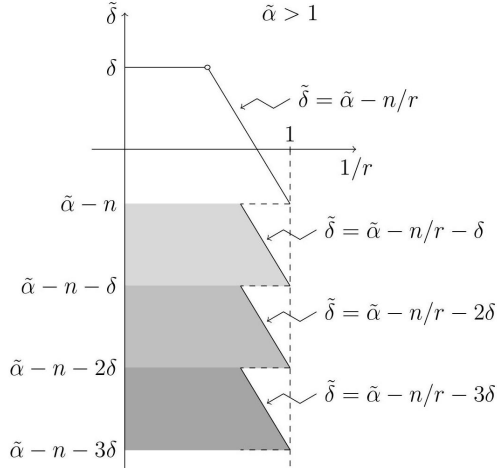
\begin{figure}[h] \centering
\begin{tikzpicture}[xscale=0.6,yscale=0.5]
\draw (0,3)--(2.5,3)--(5,-2);
\draw[dashed] (5,0)--(5,-12);
\filldraw[fill=gray!30,draw=gray!30] (0,-2)--(3.5,-2)--(5,-5)--(0,-5)--cycle;
\filldraw[fill=gray!50,draw=gray!50] (0,-5)--(3.5,-5)--(5,-8)--(0,-8)--cycle;
\filldraw[fill=gray!70,draw=gray!70] (0,-8)--(3.5,-8)--(5,-11)--(0,-11)--cycle;
\draw[dashed] (5,-2)--(3.5,-2);
\draw (3.5,-2)--(5,-5);
\draw[dashed] (5,-5)--(3.5,-5);
\draw (3.5,-5)--(5,-8);
\draw[dashed] (5,-8)--(3.5,-8);
\draw (3.5,-8)--(5,-11);
\draw[dashed] (5,-11)--(3.5,-11);
\draw[->] (-1,0)--(7,0) node[below] {$1/r$};
\draw[->] (0,-12)--(0,5) node[left] {$\tilde{\delta}$};
\draw[fill=white] (2.5,3) circle [radius=0.08];
\node[left] at (0,3) {$\delta$};
\node[left] at (0,-2) {$\tilde{\alpha}-n$};
\node[left] at (0,-5) {$\tilde{\alpha}-n-\delta$};
\node[left] at (0,-8) {$\tilde{\alpha}-n-2\delta$};
\node[left] at (0,-11) {$\tilde{\alpha}-n-3\delta$};
\draw (5,-0.1)--(5,0.1);
\node[above] at (5,0) {$1$};
\node at (5,5) {$\tilde{\alpha} > 1$};
\draw [<-] (3.5,1.5) -- ++(0.5,10pt) -- ++(0.5,-10pt) -- ++(0.5,10pt) node[right] {$\tilde{\delta} = \tilde{\alpha} - n/r$};
\draw [<-] (4.3,-3.3) -- ++(0.5,10pt) -- ++(0.5,-10pt) -- ++(0.5,10pt) node[right] {$\tilde{\delta} = \tilde{\alpha} - n/r - \delta$};
\draw [<-] (4.3,-6.3) -- ++(0.5,10pt) -- ++(0.5,-10pt) -- ++(0.5,10pt) node[right] {$\tilde{\delta} = \tilde{\alpha} - n/r - 2\delta$};
\draw [<-] (4.3,-9.3) -- ++(0.5,10pt) -- ++(0.5,-10pt) -- ++(0.5,10pt) node[right] {$\tilde{\delta} = \tilde{\alpha} - n/r - 3\delta$};
\end{tikzpicture}
\caption{Case $i)$ in the proof of Theorem \ref{teo 4.17 de tesis gladis}}
\label{FIG: pesos en los dientes}
\end{figure}

\item $\alphat-n/r-k\delta-\delta<\deltat \leq \alphat-n-k\delta$, $k \in\N \cup \{0\}$. (see figure \ref{FIG:parte de los triangulos})
\begin{figure}[h] \centering
\begin{tikzpicture}[xscale=0.6,yscale=0.5]
\draw (0,3)--(2.5,3)--(5,-2);
\filldraw[fill=gray!30,draw=gray!30] (3.5,-2)--(5,-2)--(5,-5)--cycle;
\filldraw[fill=gray!50,draw=gray!50] (3.5,-5)--(5,-5)--(5,-8)--cycle;
\filldraw[fill=gray!70,draw=gray!70] (3.5,-8)--(5,-8)--(5,-11)--cycle;
\draw (5,-2)--(3.5,-2);
\draw[dashed] (3.5,-2)--(5,-5);
\draw (5,-5)--(3.5,-5);
\draw[dashed] (3.5,-5)--(5,-8);
\draw (5,-8)--(3.5,-8);
\draw[dashed] (3.5,-8)--(5,-11);
\draw[dashed] (5,0)--(5,-12);
\draw[->] (-1,0)--(7,0) node[below] {$1/r$};
\draw[->] (0,-12)--(0,5) node[left] {$\tilde{\delta}$};
\draw[fill=white] (2.5,3) circle [radius=0.08];
\node[left] at (0,3) {$\delta$};
\node[left] at (0,-2) {$\tilde{\alpha}-n$};
\node[left] at (0,-5) {$\tilde{\alpha}-n-\delta$};
\node[left] at (0,-8) {$\tilde{\alpha}-n-2\delta$};
\node[left] at (0,-11) {$\tilde{\alpha}-n-3\delta$};
\draw (5,-0.1)--(5,0.1);
\node[above] at (5,0) {$1$};
\node at (5,5) {$\tilde{\alpha} > 1$};
\draw [<-] (3.5,1.5) -- ++(0.5,10pt) -- ++(0.5,-10pt) -- ++(0.5,10pt) node[right] {$\tilde{\delta} = \tilde{\alpha} - n/r$};
\draw [<-] (4.3,-3.3) -- ++(0.5,10pt) -- ++(0.5,-10pt) -- ++(0.5,10pt) node[right] {$\tilde{\delta} = \tilde{\alpha} - n/r - \delta$};
\draw [<-] (4.3,-6.3) -- ++(0.5,10pt) -- ++(0.5,-10pt) -- ++(0.5,10pt) node[right] {$\tilde{\delta} = \tilde{\alpha} - n/r - 2\delta$};
\draw [<-] (4.3,-9.3) -- ++(0.5,10pt) -- ++(0.5,-10pt) -- ++(0.5,10pt) node[right] {$\tilde{\delta} = \tilde{\alpha} - n/r - 3\delta$};
\end{tikzpicture}
\caption{Case $ii)$ in the proof of Theorem \ref{teo 4.17 de tesis gladis}}
\label{FIG:parte de los triangulos}
\end{figure}
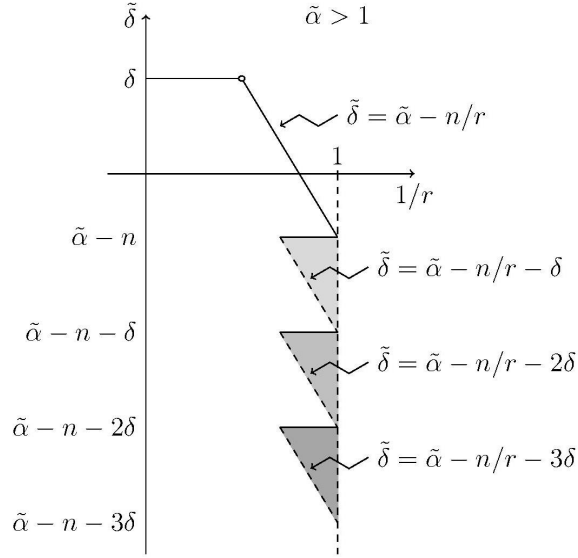
\end{enumerate}
For the case  $i)$ we consider the pairs $(w,v)$ given by $w(x)=|x|^{k\delta}$ and $v(x)=|x|^{n/r-\alphat+\deltat+k\delta}$ with
$$
\alphat-n-k\delta<\deltat\leq \min \{\alphat-n/r-k\delta,\alphat-n-k\delta+\delta\}, \;\;\;\; k\in \N.
$$
Since $v^{r'}$ satisfies the doubling condition, we use Lemma \ref{global implica local}, to estimate only global condition (\ref{parte global}).
Let $B=B(x_B,R)$ we have two cases, $|x_B|\leq R$ or $|x_B|>R$.\\
If $|x_B|\leq R$, by Proposition \ref{lema 3.5 de Gladis} and since $(n/r-\alphat+\deltat+k\delta)r'>-n$ (because $\alphat-n-k\delta<\deltat$),
$$
\int_{B}v^{r'}(x)dx \approx R^{ (n/r-\alphat+\deltat+k\delta)r'+ n}
$$
and
$$
w(B)=\int_B|x|^{k\delta}dx \approx R^{k\delta +n}.
$$
Then,
\begin{align*}
\frac{|B|^{\frac{\delta-\deltat}{n}+1}}{w(B)}&\left(\int_{\R^n\setminus B}  \frac{v^{r'}(y)}{|x_B-y|^{r'(n-\alphat+\delta)}}dy\right)^{1/r'}\\ &\menort   \frac{|B|^{\frac{\delta-\deltat}{n}+1}}{w(B)}\sum_{i=1}^{\infty}\left(\int_{2^iB\setminus 2^{i-1}B} \frac{v^{r'}(y)}{|x_B-y|^{r'(n-\alphat+\delta)}}dy\right)^{1/r'}\\
&\menort   \frac{R^{\delta-\deltat+n}}{R^{\delta k+n}}\sum_{i=1}^{\infty}\frac{(2^iR)^{n/r-\alphat+\deltat+k\delta+n/{r'}}}{(2^iR)^{n-\alphat+\delta}}\\
& \menort \sum_{i=1}^{\infty}\left(\frac{1}{2^i}\right)^{\delta-\deltat-k\delta}\\
& \approx C,
\end{align*}
where the last sum is finite because $\deltat+k \delta < \alphat-n+\delta < \delta$.

Now let $|x_B|>R$. Then there exists $N_1$ such that $\frac{|x_B|}{R}\approx 2^{N_1}$. On the other hand we have
\begin{align}
\frac{|B|^{\frac{\delta-\deltat}{n}+1}}{w(B)}&\left(\int_{\R^n\setminus B} \frac{v^{r'}(y)}{|x_B-y|^{r'(n-\alphat+\delta)}}dy\right)^{1/r'}\nonumber\\
&\menort   \frac{R^{\alphat-\deltat-n}}{|x_B|^{\delta k}}\sum_{i=1}^{\infty}\frac{1}{2^{i(n-\alphat+\delta)}}\left(\int_{2^iB}v^{r'}\right)^{1/{r'}}\label{desigualdad 5.14}
\end{align}
The last term in (\ref{desigualdad 5.14}) can be divided into $S_1$ and $S_2$ where $S_1$ is the sum up to the $N_1$-th term and $S_2$ is the sum of the remaining terms. We first estimate $S_1$
\begin{align*}
S_1&\menort  \frac{R^{\alphat-\deltat-n}}{|x_B|^{\delta k}}\sum_{i=1}^{N_1}\frac{|x_B|^{n/r-\alphat+\deltat+k\delta}(2^{i}R)^{n/{r'}}}{2^{i(n-\alphat+\delta)}}\\
&\menort \left(\frac{R}{|x_B|}\right)^{\alphat-\deltat-n/r}\sum_{i=1}^{N_1}\left(\frac{1}{2^{i}}\right)^{\delta-\deltat}\\
& \menort C,
\end{align*}
and the last sum is finite because $\deltat <\delta$.

For $S_2$ we have
\begin{align}
S_2&\menort \frac{R^{\alphat-\deltat-n}}{|x_B|^{\delta k}}\sum_{i=N_1+1}^{\infty}\frac{(2^{i}R)^{n-\alphat+\deltat+\delta k}}{2^{i(n-\alphat+\delta)}}\nonumber\\
& \approx \frac{R^{\delta k}}{|x_B|^{\delta k}}\sum_{i=N_1+1}^{\infty}\frac{1}{2^{i(\delta-\deltat-\delta k)}}.\label{desigualdad 5.15}
\end{align}
Since $\deltat +\delta k <\delta$, the last term of (\ref{desigualdad 5.15}) is less than or equal to $({R}/{|x_B|})^{\delta k}$, which is bounded by a constant.

We now estimate $ii)$. Let $\alphat -n/r-(k-1)\delta <\deltat \leq \alphat-n-k \delta$, $k\in \N\cup\{0\}$. We consider the pair $(w,v)$ defined by $w(x)=|x|^{\theta}$ and $v(x)=|x|^{\beta}$ with
$$
\theta=\alphat-n/r-k\delta-2\deltat \;\;\;\;\text{and}\;\;\;\;\beta=-k\delta-\deltat.
$$
Since $v^{r'}$ satisfies a doubling condition, by Lemma \ref{global implica local} we only  estimate the global condition (\ref{parte global}).
Let $B=B(x_B,R)$.
If $|x_B|\leq R$, by Proposition \ref{lema 3.5 de Gladis} we have
\begin{align}
\frac{|B|^{\frac{\delta-\deltat}{n}+1}}{w(B)}\bigg(\int_{\R^n\setminus B} & \frac{v^{r'}(y)}{|x_B-y|^{r'(n-\alphat+\delta)}}dy\bigg)^{1/r'}\nonumber\\ &\menort   \frac{|B|^{\frac{\delta-\deltat}{n}+1}}{w(B)}\sum_{i=1}^{\infty}\frac{1}{(2^{i}R)^{n-\alphat+\delta}}\left(\int_{2^{i}B} v^{r'}\right)^{1/r'}\label{desigualdad 5.16}\\
& \approx R^{-\deltat-\theta+\beta-n/r+\alphat}\sum_{i=1}^{\infty}\frac{1}{2^{i(n/r-\alphat+\delta-\beta)}}.\nonumber
\end{align}
Noting that
$$
-\deltat-\theta+\beta-n/r+\alphat=0\;\;\;\; \text{and}\;\;\;\; n/r-\alphat+\delta-\beta >0,
$$
it is immediate that the last sum in (\ref{desigualdad 5.16}) is bounded by a constant.

Let us now consider $|x_B|>R$. As in the case $i)$, we obtain
\begin{align*}
\frac{|B|^{\frac{\delta-\deltat}{n}+1}}{w(B)}&\left(\int_{\R^n\setminus B} \frac{v^{r'}(y)}{|x_B-y|^{r'(n-\alphat+\delta)}}dy\right)^{1/r'}\\
&\menort   \frac{R^{\alphat-\deltat}}{|x_B|^{\theta}R^{n}}\sum_{i=1}^{\infty}\frac{1}{2^{i(n-\alphat+\delta)}}\left(\int_{2^iB}v^{r'}\right)^{1/{r'}}.
\end{align*}
They, we take $S_1$ and $S_2$ as in case $i)$.\\
Since $i\leq N_1$, $\theta=\alphat -\deltat+\beta-n/r$ and $|x_B|>2^{i}R$, we have
\begin{align*}
S_1&\menort \frac{R^{\alphat-\deltat}}{R^n}\sum_{i=0}^{N_1}\frac{|x_B|^{\beta-\theta}(2^iR)^{n/{r'}}}{2^{i(n-\alphat+\delta)}}\\
&\menort R^{\alphat-\deltat-n/r-\theta+\beta}\sum_{i=0}^{N_1}2^{-i(n-\alphat+\delta-n/{r'}+\theta-\beta)}\\
&\menort \sum_{i=1}^{N_1} 2^{i(\deltat-\delta)},
\end{align*}
which is finite since $\deltat<\delta$.

For $S_2$ we get
\begin{align*}
S_2&\menort \frac{R^{\alphat-\deltat+n/{r'}}}{|x_B|^{\theta}R^n}\sum_{i=N_1+1}^{\infty}\frac{|x_B|^{\beta}}{2^{i(n/r-\alphat+\delta)}}\\
&\menort \frac{R^{\alphat-\deltat+\beta-n/r}}{|x_B|^{\theta}}\sum_{i=N_1+1}^{\infty}2^{-i(n/r-\alphat+\delta-\beta)}\\
&\menort \left(\frac{R}{|x_B|}\right)^{\theta}\sum_{i=N_1+1}^{\infty}2^{-i(n/r-\alphat+\delta-\beta)}\\
\end{align*}
Now, since $|x_B|>R$, $\alphat-\deltat+\beta-n/r=\theta>0$ and $n/r-\alphat+\delta-\beta >0$, we obtain
$$
S_2\menort C.
$$
 This concludes the proof of $ii)$.

For the case $r=1$ and $\deltat< \alphat-n$ we set $w(x)=|x|^{-\deltat}$
 and $v(x)=|x|^{n-\alphat}$.
By Lemma \ref{global implica local}, we shall estimate (\ref{parte global}). Let $B=B(x_B,R)$, if $|x_B|\leq R$, we then
\begin{align*}
\frac{|B|^{(\delta-\deltat)/n+1}}{w(B)}\left\|\frac{\chi_{\R^n\setminus B} v}{(|B|^{1/n}+|x_B-.|)^{n-\alphat+\delta}}\right\|_{\infty}&\menort \frac{|B|^{(\delta-\deltat)/n+1}}{w(B)}\sum_{i=1}^{\infty}\frac{1}{(2^{i}R)^{n-\alphat+\delta}}\|\chi_{B_i} v\|_{\infty}\\
&\menort \frac{R^{\delta-\deltat+n}}{R^{n-\deltat}}\sum_{i=1}^{\infty}\frac{(2^{i}R)^{n-\alphat}}{(2^{i}R)^{n-\alphat+\delta}}\\
&\approx C \sum_{i=1}^{\infty} \frac{1}{2^{i\delta}}\\
&\menort C.
\end{align*}
If $|x_B|>R$, we proceed as in the case $p>1$ to obtain that the first term of the above inequality is bounded by $S_1$ and $S_2$ where
$$
S_1\approx \frac{R^{\delta-\deltat+n}}{|x_B|^{-\deltat}R^{n}}\sum_{i=1}^{N_1}\frac{\|\chi_{B_i}v\|_{\infty}}{(2^{i}R)^{n-\alphat+\delta}},
$$
$$
S_2\approx \frac{R^{\delta-\deltat+n}}{|x_B|^{-\deltat}R^{n}}\sum_{i=N_1+1}^{\infty}\frac{\|\chi_{B_i}v\|_{\infty}}{(2^{i}R)^{n-\alphat+\delta}}.
$$
In order to estimate $S_1$, since $|x_B|>2^{i}R$ for $i\leq N_1$, we have
$$
S_1\menort \frac{R^{\delta-\deltat}}{|x_B|^{-\deltat}}\sum_{i=1}^{N_1}\frac{|x_B|^{n-\alphat}}{(2^{i}R)^{n-\alphat+\delta}}\menort R^{\alphat-\deltat-n}|x_B|^{\deltat+n-\alphat}\menort C.
$$
On the other hand

$$
S_2\menort \frac{R^{\delta-\deltat}}{|x_B|^{-\deltat}}\sum_{i=N_1+1}^{\infty}\frac{(2^{i}R)^{n-\alphat}}{(2^{i}R)^{n-\alphat+\delta}}\menort \left(\frac{R}{|x_B|}\right)^{-\deltat}\sum_{i=2}^{\infty}\frac{1}{2^{i\delta}},
$$
and since $\deltat<\alphat-n<0$ and $|x_B|>R$, the last term is bounded by a constant.
\end{dem}

\begin{prop} \label{clase no abiertas} Let $0\leq \alpha<n$, $0<\delta<\min\{\eta, (n-\alpha)/m\}$. Let  $1\leq r <\infty$ and $\deltat\leq \min\{\delta, \alphat-n/r\}$, excluding the case $\deltat=\delta$ when $\alphat-n/r=\delta$. Then there exist pairs of weights $(w,v)$ belonging to $\H(r,\alphat,\deltat)$ such that $(w,v)$ does not belong to $\H((r't)',\alphat,\deltat)$ for any $t>0$, with $t\neq 1$.
\end{prop}
\begin{dem}
We shall exhibit  examples of pairs of weights $(w,v)$ such that $(w,v) \in \H(r,\alphat,\deltat)$ but $(w,v) \notin \H((r't)',\alphat,\deltat)$ for all $t>0$, $t\neq 1$.\\
Let  $1< r <\infty$ and $\alphat-n< \deltat \leq \alphat-n/r< \delta$. Let us consider the pair of weights $(w,v)=(1,|x|^{-\theta})$ with $\theta=\alphat-n/r-\deltat$. By Theorem 3.6 of \cite{PraReinf} we have that  $(w,v) \in \mathbb{H}(r,\alphat,\deltat)\subset \H(r,\alphat,\deltat)$.
Let us see that $(w,v)$ does not belong to $\H((r't)',\alphat,\deltat)$ for any $t>0$, with $t\neq1$.
By Lemma  \ref{lema condicion de peso es equivalente a la parte local mas la global}, it is enough to show that the pair $(w,v)$ does not satisfy the local condition
\begin{equation}\label{parte local especial}
|B|^{\frac{\alphat-\deltat}{n}-1}\|\chi_B v\|_{r't}\menort \frac{w(B)}{|B|}.
\end{equation}
Let $0<t<n/(\theta r')$ (otherwise $\|\chi_B v\|_{r't}=\infty$). Let $B=B(0,R)$, then the left hand side of the inequality (\ref{parte local especial}) is bounded below by

$$
|B|^{-1+(\alphat-\deltat)/n}\|v\chi_B\|_{r't}\mayort R^{-n+\alphat-\deltat}R^{-\theta+n/(r't)}\mayort R^{-n/(r't')}
$$
and the last term tends  to infinity when $R$ tends to zero or infinity if $t>1$ or $t<1$, respectively. \\
We now consider the case $1< r <\infty$ and $\deltat\le\delta<\alphat-n/r$, and the pair $(|x|^{-\beta},|x|^{-\theta})$, where
$$
\alphat-n/r-\delta<\theta<n/{r'} \;\;\;\;\;\;\text{and}\;\;\;\;\;\; 0<\beta=\theta +\delta-\alphat+n/r.
$$
By Theorem 3.6 of \cite{PraReinf} we have that  $(w,v) \in \mathbb{H}(r,\alphat,\deltat)\subset \H(r,\alphat,\deltat)$. Let us now see that it does not belong to $\H((r't')',\alphat,\delta)$ for any $t>0$, with $t\neq1$.   By Lemma \ref{lema condicion de peso es equivalente a la parte local mas la global} it is enough to see that there exists a ball $B$ such that the local condition (\ref{parte local especial}) does not hold. In fact, if $B=B(0,R)$ then
$$
\frac{|B|^{(\alphat-\deltat)/n}\|v\chi_B\|_{r't}}{w(B)}\mayort R^{\beta-n+\alphat-\delta}R^{-\theta+n/(r't)}\simeq R^{-n/{r'}+n/(r't)}\simeq R^{-n/(r't')}.
$$
Consequently, the last term tends to infinity when $R$ tends to zero or infinity if $t>1$ or $t<1$, respectively.\\
Let $1<r<\infty$, $\alphat-n-k\delta<\deltat\leq \min \{\alphat-n/r-k\delta, \alphat-n-(k-1)\delta\}$, $k\in \N$ and $(w,v)=(|x|^{k\delta}, |x|^{\theta})$, where
$$\theta =n/r-\alphat+\deltat+k\delta.$$
By Theorem \ref{teo 4.17 de tesis gladis}, we have $(w,v) \in \H(r,\alphat,\deltat)$. However, $(w,v)\notin \H((r't')',\alphat,\delta)$ for any $t>0$, with $t\neq1$ because the local condition (\ref{parte local especial}) does not hold. In fact, if $B=B(0,R)$,
$$
\frac{|B|^{(\alphat-\deltat)/n}\|v\chi_B\|_{r't}}{w(B)}\mayort R^{\alphat-\deltat-k\delta-n}R^{n/r-\alphat+\deltat+k\delta+n/{r't}}\simeq R^{-n/{r't'}}
$$
where the last term tends to infinity when $R$ tends to zero or infinity if $t>1$ or $t<1$, respectively.\\

Similar arguments show that, if $1<r<\infty$ and $\alphat-n/r-k\delta-\delta<\deltat\leq \alphat-n-k\delta$, $k\in \N$, then the pair  $(w,v)=(|x|^{\theta}, |x|^{\beta})$, where
$$
\theta=\alphat-n/r-k\delta-2\deltat \;\;\;\;\text{and}\;\;\;\;\beta=-k\delta-\deltat,
$$
belongs to $\H(r,\alphat,\deltat)$. However, $(w,v)\notin \H((r't')',\alphat,\delta)$ for any $t>0$, with $t\neq1$ because the local condition (\ref{parte local especial}) does not hold.\\

For the case $r=1$ and $\deltat=\alphat-n$, it is immediate that the pair $(w,v)$ given by $w=v=1$ belongs to $\H(1,\alphat,\alphat-n)$. However, it is easy to check that $(w,v) \notin \H(1+\epsilon,\alphat,\alphat-n)$ for every $\epsilon>0$.

Finally, if $r=1$ and $\deltat< \alphat-n$, let us consider the pair $(|x|^{-\deltat},|x|^{n-\alphat})$. It was proved in Theorem \ref{teo 4.17 de tesis gladis} that $(w,v)$ belongs to $\H(1,\alphat, \deltat)$. Let us see that $(w,v)$ does not belong to $\H(1+\epsilon, \alphat, \deltat)$ for any $\epsilon >0$. By Lemma \ref{lema condicion de peso es equivalente a la parte local mas la global} it is enough to show that $(w,v)$ does not satisfy condition (\ref{parte local}) with $r=1+\epsilon$. In fact, if $B=B(0,R)$, we get
$$
\frac{|B|^{(\alphat-\deltat)/n}}{w(B)}\|v\chi_B\|_{(1+\epsilon)'}\mayort R^{n/(1+\epsilon)'}
$$
and the last expression tends to $\infty$ when $R$ tends to $\infty$.
\end{dem}

\section{Proof of the main results}
We now give some previous lemmas that we shall  use in the proofs of the main results. We are considering $m\in \N\cup \{0\}$.
\begin{lema}\label{acotacion de la I3}
Let  $0\leq\alpha<n$, $0<\delta<\min(\eta, (n-\alpha)/m)$ and $1\leq r \leq \infty$. Let $K_\alpha\in  K^{*}_{\alpha, \infty}$ and $b\in \Lambda(\delta)$. If $(w,v)\in \H(r,\alphat,\deltat)$, then
$$
\int_{(2B)^c} |b(x)-b(z)|^m |K_\alpha(x-z)-K_\alpha(y-z)||f(z)| dz \,\menort\, \|b\|^m_{\Lambda(\delta)}w(B)|B|^{\deltat/n-1}\left\|f/v\right\|_{r}
$$
for all $x,y \in B$.
\end{lema}

\begin{dem}
If $x,y \in B$, by using that  $b\in\Lambda(\delta)$ and $K_\alpha\in K^{*}_{\alpha, \infty}$, we have that
\begin{align*}
\int_{(2B)^c} & |b(x)-b(z)|^m |K_\alpha(x-z)-K_\alpha(y-z)||f(z)| dz\\
&\menort \|b\|^m_{\Lambda(\delta)}\int_{(2B)^c} |x-z|^{\delta m} |K_\alpha(x-z)-K_\alpha(y-z)||f(z)| dz\\
&\menort \|b\|^m_{\Lambda(\delta)}\sum_{j=1}^{\infty}\frac{2^{j\delta m}|B|^{\delta m/n+\eta/n}}{2^{j(n-\alpha+\eta)}|B|^{(n-\alpha+\eta)/n}}\int_{2^{j+1}B\setminus 2^{j}B }|f(z)|dz\\
&\menort \|b\|^m_{\Lambda(\delta)}|B|^{\frac{\delta m-n+\alpha}{n}}\sum_{j=1}^{\infty}(2^{j})^{\delta m -n+\alpha-\eta}\int_{2^{j+1}B\setminus 2^{j}B }|f(z)|vv^{-1}dz.
\end{align*}
Now, we can apply  H\"older's inequality to get
\begin{align*}
\int_{(2B)^c} &|b(x)-b(z)|^m |K_\alpha(x-z)-K_\alpha(y-z)||f(z)| dz\\
&\menort \|b\|^m_{\Lambda(\delta)}\left\|f/v\right\|_{r}|B|^{\frac{\delta m-n+\alpha}{n}}\sum_{j=1}^{\infty}2^{j(\delta m-n+\alpha-\eta)}\left(\int_{2^{j+1}B\setminus 2^{j}B }v^{r'}(z)dz\right)^{1/{r'}}\\
&\menort \|b\|^m_{\Lambda(\delta)}\left\|f/v\right\|_{r}|B|^{\frac{\delta m-n+\alpha}{n}}\sum_{j=1}^{\infty}2^{j(\delta m-n+\alpha-\eta)}|2^jB|^{\frac{n-\alphat+\delta}{n}}\times\\
&\;\;\;\;\;\;\;\;\;\;\;\;\;\;\;\;\;\;\;\;\;\;\;\;\;\;\;\;\;\;\;\;\;\;\;\;\;\;\;\;\;\;\;\;\;\;\;\;\times\left(\int_{2^{j+1}B\setminus 2^{j}B }\frac{v^{r'}(z)}{|x_B-z|^{r'(n-\alphat+\delta)}}dz\right)^{1/{r'}}\\
&\menort \|b\|^m_{\Lambda(\delta)}\left\|f/v\right\|_{r}|B|^{\frac{\delta}{n}}\sum_{j=1}^{\infty}2^{j(\delta-\eta)}\left(\int_{2^{j+1}B\setminus 2^{j}B }\frac{v^{r'}(z)}{|x_B-z|^{r'(n-\alphat+\delta)}}dz\right)^{1/{r'}}\\
&\menort \|b\|^m_{\Lambda(\delta)}\left\|f/v\right\|_{r}|B|^{\frac{\delta}{n}}\sum_{j=1}^{\infty}2^{j(\delta-\eta)}\left(\int_{\R^n\setminus B }\frac{v^{r'}(z)}{|x_B-z|^{r'(n-\alphat+\delta)}}dz\right)^{1/{r'}}
\end{align*}
Since $\sum_{j=1}^{\infty}(2^{j})^{\delta-\eta}$ is finite, by using the global condition (\ref{parte global}), we get
$$
\int_{(2B)^c} |b(x)-b(z)|^m |K_\alpha(x-z)-K_\alpha(y-z)||f(z)| dz
\menort \|b\|^m_{\Lambda(\delta)}|B|^{\frac{\deltat}{n}}\left\|f/v\right\|_{r}\frac{w(B)}{|B|}.
$$
\end{dem}

\begin{lema}\label{acotacion de la I1}
Let $0<\alpha<n$,  $0<\delta<\min(\eta, (n-\alpha)/m)$ and $1\leq r \leq \infty$. Let $K_\alpha\in S_{\alpha}^*$ and $b\in \Lambda(\delta)$. Let $(w,v)\in \H(r,\alphat, \deltat)$, then
$$
\frac{1}{w(B)}\int_{B}  |T^m_{\alpha,b}f\chi_{2B}(x)| dx \,\menort \|b\|^m_{\Lambda(\delta)}|B|^{\deltat/n}\left\|f/v\right\|_{L^r}
$$
for all ball $B  \in \R^n$.
\end{lema}

\begin{dem}

By using Tonelli's theorem and the fact that $b\in \Lambda(\delta)$, we obtain that
\begin{align*}
\frac{1}{w(B)}\int_{B}  |T^m_{\alpha,b}&f\chi_{2B}(x)| dx \\
\leq&\frac{1}{w(B)}\int_{B}\int_{2B} |b(x)-b(y)|^m |K_\alpha(x-y)||f(y)| dy\,dx\\
&\menort \frac{1}{w(B)}\int_{2B}|f(y)|\left(\int_{B} |b(x)-b(y)|^m |K_\alpha(x-y)| dx\right)dy\\
&\menort \|b\|_{\Lambda(\delta)}^m \frac{1}{w(B)}\int_{2B}|f(y)|\left(\int_{2B}|x-y|^{\delta m} |K_\alpha(x-y)| dx\right)dy.
\end{align*}
Now, since $K_\alpha \in S_{\alpha}^*$ and $\delta m < n-\alpha$, it is easy to see that

$$
\int_{2B}|x-y|^{\delta m} |K_\alpha(x-y)| dx\menort |B|^{(\delta m+\alpha)/n}.
$$

Then, by H\"older's inequality, using Lemma \ref{lema del doble v}, we have

\begin{align*}
\frac{1}{w(B)}\int_{B}  |T^m_{\alpha,b}f\chi_{2B}(x)| dx
& \menort \|b\|_{\Lambda(\delta)}^m|B|^{\frac{\delta m+\alpha}{n}}\|v\chi_{2B}\|_{r'}\left\|f/v\right\|_{r}(w(B))^{-1}\\
& \menort \|b\|_{\Lambda(\delta)}^m|B|^{\frac{\delta m+\alpha}{n}}\left\|f/v\right\|_{r}|B|^{\frac{\deltat -\alphat}{n}}\\
& \menort \|b\|_{\Lambda(\delta)}^m|B|^{\frac{\deltat}{n}}\left\|f/v\right\|_{r}.\\
\end{align*}
\end{dem}

 It was proved in \cite{LuchiGla} that if  $0<\delta<\min(\eta, n/m)$, $1< p < n/m$, $\frac{1}{q}=\frac{1}{p}-\frac{m\delta}{n}$, $b\in \Lambda(\delta)$, then
\begin{equation}\label{acotacion sin pesos para integrales singulares}
\left ( \int_{\R^n}|T^m_b f(x)|^q dx\right)^{1/q} \,\menort \|b\|_{\Lambda(\delta)} \left ( \int_{\R^n}| f(x)|^p dx\right)^{1/p}
\end{equation}
The result above will be useful in the following estimate.
\begin{lema}\label{acotacion de la I1 para alpha cero}
Let $0<\delta<\min(\eta, n/m)$ and $1\leq r \leq \infty$. Let $K_0\in S_{0}^*$ and $b\in \Lambda(\delta)$. Let $(w,v)\in \H(r,m\delta, \deltat)$ and $v^{r'}\in RH(\frac{p(r-1)}{r-p})$ for some $p\in(1,r)$, then
$$
\frac{1}{w(B)}\int_{B}  |T^m_{b}f\chi_{2B}(x)| dx \,\menort \|b\|^m_{\Lambda(\delta)}|B|^{\deltat/n}\left\|f/v\right\|_{L^r}
$$
for all ball $B  \in \R^n$.

\end{lema}
\begin{dem}

$$
\frac{1}{w(B)}\int_{B}  |T^m_{b}f\chi_{2B}(x)| dx\;= \;\frac{|B|}{w(B)}\left(\frac{1}{|B|}\int_{B}  |T^m_{b}f\chi_{2B}(x)| dx\right)
$$
Let $s \in (p,r)$ such that $\frac{1}{s}=\frac{1}{p}-\frac{\delta m}{n}$, by Jessen and (\ref{acotacion sin pesos para integrales singulares}), we have

\begin{align*}
\frac{1}{w(B)}\int_{B}  |T^m_{b}f\chi_{2B}(x)| dx & \menort \frac{|B|^{1-1/s}}{w(B)}\left(\int_{\R^n}  |T^m_{b}f\chi_{2B}(x)|^s dx\right)^{1/s} \\
&\menort \frac{|B|^{1-1/s}}{w(B)}\|b\|^m_{\Lambda(\delta)}\left(\int_{2B}  |f(x)|^p dx\right)^{1/p}\\
& \menort \frac{|B|^{1-1/s+1/p}}{w(B)}\|b\|^m_{\Lambda(\delta)}\left(\frac{1}{|2B|}\int_{2B}  \frac{|f(x)|^p v^p}{v^p} dx\right)^{1/p},\\
\end{align*}
Now, first applying Hölder $\left( {r/p},\left({r/p}\right)^{'}\right)$  then using the fact that  $v^{r'}\in RH(\frac{p(r-1)}{r-p})$ and finally using the Lemma \ref{lema del doble v}, we obtain

\begin{align*}
\frac{1}{w(B)}\int_{B}  |T^m_{b}f\chi_{2B}(x)| dx &  \menort \frac{|B|^{1+\frac{\delta m}{n}}}{w(B)}\|b\|^m_{\Lambda(\delta)}\left(\frac{1}{|2B|}\int_{2B}  \frac{|f(x)|^r}{v^r} dx\right)^{\frac{1}{r}}\left(\frac{1}{|2B|}\int_{2B}  v^{p(\frac{r}{p})^{'}} dx\right)^{\frac{1}{p(\frac{r}{p})^{'}}}\\
& \menort \frac{|B|^{1+\frac{\delta m}{n}-\frac{1}{r}}}{w(B)}\|b\|^m_{\Lambda(\delta)} \|f/v\|_{L^r}\left(\frac{1}{|2B|}\int_{2B}  v^{r'} dx\right)^{1/{r'}}\\
& \menort |B|^{\deltat/n}\|b\|^m_{\Lambda(\delta)} \|f/v\|_{L^r} .\\
\end{align*}

\end{dem}

We now proceed with the proof of Theorems \ref{teo para integrales singulares} and \ref{teorema A}.

\medskip

\begin{prueba}{Theorems \ref{teo para integrales singulares} and \ref{teorema A}}
Let $f/v \in L^r(\R^n)$. Let $B\subset\R^n$ be a ball and $x\in B$. We split  $f=f_1+f_2$ with $f_1=f\chi_{2B}$ and define $a_B= \frac{1}{|B|}\int_{B} T^{m}_{\alpha,b}f_2$. Then,
\begin{align*}
\frac{1}{w(B)}\int_B|T^m_{\alpha,b}f(x)-a_B|\,dx &\menort \frac{1}{w(B)}\int_B|T^m_{\alpha,b}f_1(x)|\,dx + \frac{1}{w(B)}\int_B|T^m_{\alpha,b}f_2(x)-a_B|\,dx \\
&= I_1+I_2
\end{align*}
By Lemma \ref{acotacion de la I1} if $0<\alpha<n$ or by Lemma \ref{acotacion de la I1 para alpha cero} if $\alpha=0$, we have
\begin{equation}
I_1\menort \|b\|^m_{\Lambda(\delta)}|B|^{\deltat/n}\left\|f/v\right\|_{r}.
\end{equation}
 Since
$$
|T^m_{\alpha,b}f_2(x)-a_B|= |T^m_{\alpha,b}f_2(x)-(T^m_{\alpha,b}f_2)_B|\menort \frac{1}{|B|}\int_B |T^m_{\alpha,b}f_2(x)-T^m_{\alpha,b}f_2(y)|dy,
$$
then,
\begin{equation}
I_2\menort \frac{1}{w(B)}\int_B\frac{1}{|B|}\int_B |T^m_{\alpha,b}f_2(x)-T^m_{\alpha,b}f_2(y)|dy\;dx.
\end{equation}
Let $A=|T^m_{\alpha,b}f_2(x)-T^m_{\alpha,b}f_2(y)|$. If $x,y \in B$
\begin{align*}
A &\menort \int_{(2B)^c}|(b(x)-b(z))^m K_\alpha (x-z)-(b(y)-b(z))^mK_\alpha (y-z)||f(z)|dz \\
&\menort \int_{(2B)^c}|b(x)-b(z)|^m |K_\alpha (x-z)-K_\alpha (y-z)||f(z)|dz\\
&+\int_{(2B)^c}|(b(x)-b(z))^m-(b(y)-b(z))^m ||K_\alpha (y-z)||f(z)|dz\\
&= I_3+I_4.
\end{align*}
By Lemma \ref{acotacion de la I3}, we have
$$
I_3\menort \|b\|^m_{\Lambda(\delta)}w(B)|B|^{\frac{\deltat}{n}-1}\left\|f/v\right\|_{r}.
$$
In order to estimate $I_4$, we use that $b\in\Lambda(\delta)$. If $A_j=2^{j+1}B\setminus 2^jB $, then
\begin{align*}
I_4 &\menort |b(x)-b(y)|\sum_{k=0}^{m-1}\int_{(2B)^c}|b(x)-b(z)|^{m-1-k}|b(y)-b(z)|^k|K_\alpha(y-z)||f(z)|dz \\
& \menort \|b\|_{\Lambda(\delta)}^m |B|^{\delta/n}\sum_{j=1}^{\infty}|2^{j+1}B|^{{\delta(m-1)}/n}\int_{A_j}|K_\alpha(x-z)||f(z)|dz\\
& \menort \|b\|_{\Lambda(\delta)}^m |B|^{\delta/n}\sum_{j=1}^{\infty}\int_{A_j}\frac{|f(z)|}{|x_B-z|^{(n-\alphat+\delta)}}dz\\
& \menort \|b\|_{\Lambda(\delta)}^m |B|^{\delta/n}\int_{\R^n\setminus B}\frac{|f(z)|}{|x_B-z|^{(n-\alphat+\delta)}}dz\\
\end{align*}
Then, by H\"older's inequality, and global condition (\ref{parte global}), we have

\begin{align*}
I_4 & \menort \|b\|_{\Lambda(\delta)}^m |B|^{\delta/n}\|f/v\|_{r}\left(\int_{\R^n\setminus B}\frac{v^{r'}(z)}{|x_B-z|^{r'(n-\alphat+\delta)}}dz\right)^{1/{r'}}\\
& \menort \|b\|_{\Lambda(\delta)}^m |B|^{\deltat/n}\|f/v\|_{r}\frac{w(B)}{|B|}
\end{align*}
We now get
\begin{align*}
I_2 & \menort \frac{1}{w(B)} \int_B \frac{1}{|B|}\int_B (I_3+I_4)dydx\\
&\menort \|b\|^m_{\Lambda(\delta)}|B|^{\frac{\deltat}{n}}\left\|f/v\right\|_{r}\\
\end{align*}
%
%
%
%
 and then, we obtain that
$$
\frac{1}{w(B)}\int_B|T^m_{\alpha,b}f(x)-a_B|\,dx\menort \|b\|^m_{\Lambda(\delta)}|B|^{\frac{\deltat}{n}}\left\|f/v\right\|_{r},
$$
so it remains to take supremum over all the balls $B$ to get the desired result.
\end{prueba}

\small
\markright{}

\bibliographystyle{abbrv}
\bibliography{PRADO}
\noindent Gladis Pradolini, {\sl CONICET and Departamento de Matem\'atica, Facultad de Ingenier\'ia Qu\'imica, UNL, 3000, Santa Fe, Argentina}.\\
\noindent e-mail address: gladis.pradolini@gmail.com

\noindent Jorgelina Recchi, {\sl Departamento de Matem\'aticas, Universidad Nacional del Sur (UNS), 8000, Bah\'ia Blanca, Argentina.}\\
\noindent e-mail address: drecchi@uns.edu.ar, jrecchi@gmail.com
\end{document}